\pdfoutput=1
\RequirePackage{ifpdf}
\ifpdf % We~are running pdfTeX in pdf mode
\documentclass[pdftex]{sigma}
\else
\documentclass{sigma}
\fi

\numberwithin{equation}{section}

\newtheorem*{Theorem*}{Theorem}

 { \theoremstyle{definition}

 }

\newmuskip\pFqskip
\mathchardef\pFcomma=\mathcode`, % keep a copy of the comma
\newcommand*\pFq[5]{%
 \begingroup
 \begingroup\lccode`~=`,
 \lowercase{\endgroup\def~}{\pFcomma\mkern\pFqskip}%
 \mathcode`,=\string"8000
 {}_{#1}F_{#2}\biggl(\genfrac..{0pt}{}{#3}{#4};#5\biggr)%
 \endgroup
}

\newcommand{\expp}[1]{\text{exp}\left(#1\right)}

\newcommand{\pp}[1]{\left({#1}\right)}
\newcommand{\bb}[1]{\left[{#1}\right]}
\usepackage{tikz}
\usepackage{tikz-cd}

\begin{document}
\allowdisplaybreaks

\newcommand{\arXivNumber}{2209.10725}

\renewcommand{\PaperNumber}{090}

\FirstPageHeading

\ShortArticleName{Para-Bannai--Ito Polynomials}

\ArticleName{Para-Bannai--Ito Polynomials}

\Author{Jonathan PELLETIER~$^{\rm a}$, Luc VINET~$^{\rm ab}$ and Alexei ZHEDANOV~$^{\rm c}$}

\AuthorNameForHeading{J.~Pelletier, L.~Vinet and A.~Zhedanov}

\Address{$^{\rm a)}$~Centre de Recherches Math\'ematiques, Universit\'e de Montr\'eal,\\
\hphantom{$^{\rm a)}$}~P.O. Box 6128, Centre-ville Station,
Montr\'eal (Qu\'ebec), H3C 3J7, Canada}
\EmailD{\href{mailto:jonathanpelletier9179@gmail.com}{jonathanpelletier9179@gmail.com}, \href{mailto:luc.vinet@umontreal.ca}{luc.vinet@umontreal.ca}}

\Address{$^{\rm b)}$~IVADO, Montr\'eal (Qu\'ebec), H2S 3H1, Canada}

\Address{$^{\rm c)}$~School of Mathematics, Renmin University of China, Beijing 100872, P.R. China}
\EmailD{\href{mailto:zhedanov@yahoo.com}{zhedanov@yahoo.com}}

\ArticleDates{Received June 10, 2023, in final form October 28, 2023; Published online November 10, 2023}

\Abstract{New bispectral polynomials orthogonal on a Bannai--Ito bi-lattice (uniform quadri-lattice) are obtained from an unconventional truncation of the untruncated Bannai--Ito and complementary Bannai--Ito polynomials. A complete characterization of the resulting para-Bannai--Ito polynomials is provided, including a three term recurrence relation, a Dunkl-difference equation, an explicit expression in terms of hypergeometric series and an orthogonality relation. They are also derived as a $q\to -1$ limit of the $q$-para-Racah polynomials. A connection to the dual $-1$ Hahn polynomials is also established.}


\Classification{33C45}

\section{Introduction}
Following the classification of Leonard \cite{Leo_class}, Bannai and Ito introduced a new classification of orthogonal polynomials as the polynomials that satisfy the Leonard duality property \cite{BIbook}. With it, they introduced the first $-1$ orthogonal polynomials: the Bannai--Ito polynomials. Vinet, Zhedanov and collaborators have since lead the way in research on $-1$ orthogonal polynomials starting with \cite{L-1J,B-1J} and giving an almost complete description of the $-1$ part of the Askey scheme. In parallel, the same group of people introduced a new category of orthogonal polynomials: the para-polynomials~\cite{pararacah,para_kraw}. The two categories are, a priori, not exclusive, which raise the question of the existence of a family of orthogonal polynomials at the crossing of the two. In this paper, we answer that question at the affirmative by presenting the para-Bannai--Ito polynomials.

There are many reasons to be interested in a $q\to -1$ limit of the $q$-para-Racah polynomials. From a mathematical physics standpoint, since the discovery of the para-Krawtchouk polynomials \cite{para_kraw}, para-polynomials have been linked to perfect state transfer (PST) and fractional revival (FR) in $XX$-spin chains with nearest neighbour couplings (see also \cite{frac_rev_pararacah}). A necessary feature for PST is the persymmetry (i.e., symmetry under anti-diagonal reflections) of the underlying Jacobi matrix and, for FR \cite{persymmetric}, the existence of an isospectral deformation of this persymmetric matrix. The models with these state transport properties that are based on $-1$ orthogonal polynomials have not been developed. From a mathematical viewpoint, para-polynomials have appeared as basis functions for the finite-dimensional representations of algebras of the Sklyanin type \cite{sky_qpK,rat_sky,skly_qpara}. The $q$-para-Racah polynomials arise in particular in connection with the degenerate Sklyanin algebra \cite{Gorsky_1993}. The $q\to -1$ limits of these $q$-para-Racah polynomials therefore stand to provide a representation basis for a $-1$ version of this algebra. The study of the para-Bannai--Ito polynomials is also of interest in the elaboration of the scheme of -1 orthogonal polynomials, but also to obtain the complete set of para-polynomials; Indeed, of the three types of grid obtained from the general Askey--Wilson grid in the cases $|q|\neq 1$, $q=1$ and $q=-1$ , only the Bannai--Ito grid ($q=-1$) has not been investigated as a potential candidate to yield a new family of orthogonal polynomials of the para type (appreciating that the quadratic bi-lattice grid leads to the para-Racah polynomials \cite{pararacah} and that the $q$-quadratic bi-lattice grid takes one to the $q$-para-Racah polynomials \cite{qpararacah}). As such, the category of para-polynomials is characterized by three main properties. Each family of para-polynomials is a finite system of orthogonal polynomials obtained from a unconventional truncation of a classical family, the Jacobi matrix representing the recurrence relation is persymmetric or an isospectral deformation of a persymmetric matrix and the polynomials are orthogonal on the union of two regular lattices (regular lattices are limits of the $q$-quadratic lattice).

The goal of this paper is to fully characterize the $q\to -1$ limit of the $q$-para-Racah polynomials. We name them para-Bannai--Ito polynomials, as they are obtained from an unconventional truncation of the general Bannai--Ito \cite{BI} and general complementary Bannai--Ito polynomials~\cite{CBI}. The paper has the following structure. In Section \ref{C2_GCBI_GBI}, we review some important properties of the general Bannai--Ito and complementary Bannai--Ito polynomials, as well as their relations via Christoffel/Geronimus transformations \cite{haine,zhedanov_rational_1997}. In Section \ref{sec3}, the truncation condition for the general complementary Bannai--Ito with $N=2j$ ($j$ even) is presented with the corresponding three-term recurrence relation, the Dunkl-difference equation, the explicit expression in terms of hypergeometric series and the orthogonality relation. A Geronimus transformation is used to obtain, in Section \ref{sec4}, the recurrence relation of the para-Bannai--Ito polynomials for $N=2j+1$ ($j$ even) through the same truncation, albeit in the case of the general Bannai--Ito polynomials. The corresponding cases for $j$ odd are treated in Appendixes~\ref{sec:appA} and~\ref{sec:appB}. In Section \ref{sec5}, the para-Bannai--Ito polynomials are obtained as a $q\to -1$ limit of the $q$-para-Racah polynomials and a~connection to the dual $-1$ Hahn polynomials is presented.

\section{General Complementary Bannai--Ito and Bannai--Ito}\label{C2_GCBI_GBI}
Let us first review some properties of the general untruncated complementary Bannai--Ito polynomials. Together with the general Bannai--Ito polynomials, they sit atop the $q\to -1$ limit of the $q$-Askey scheme and depend on four parameters $\rho_1$, $\rho_2$, $r_1$, $r_2$. They admit an explicit expression given by
\begin{gather}
 \mathbf{I}_{2n}(x;\rho_1,\rho_2,r_1,r_2) = l_n^{(1)}\mathbf{W}_n\pp{({\rm i}x)^2;\rho_2,\rho_1+1,-r_1+\frac{1}{2},-r_2+\frac{1}{2}}\label{C2_I2n},\\
 \mathbf{I}_{2n+1}(x;\rho_1,\rho_2,r_1,r_2) = (x-\rho_2)l_n^{(2)}\mathbf{W}_n\pp{({\rm i}x)^2;\rho_2+1,\rho_1+1,-r_1+\frac{1}{2},-r_2+\frac{1}{2}}\label{C2_I2n1},
\end{gather}
where $l_n^{(i)}$ are normalizing factor to ensure monicity, $\mathbf{W}_n\big(x^2\big)$ are Wilson polynomials given by
\begin{align}
 \mathbf{W}_n\bigl(x^2;a,b,c,d\bigr) &= \pFq{4}{3}{-n,n+a+b+c+d-1,a-{\rm i}x,a+{\rm i}x}{a+b,a+c,a+d}{1}\nonumber\\
 &= \sum_{k=0}^{n}\frac{(-n,n+a+b+c+d-1,a-{\rm i}x,a+{\rm i}x)_{k}}{(1)_{k}(a+b,a+c,a+d)_{k}},\label{decompositionwilson}
\end{align}
and $(a)_{k} = (a)(a+1)\cdots(a+k-1)$ and $(a_1,a_2,\dots,a_n)_{k} = (a_1)_{k}(a_2)_{k}\cdots(a_n)_{k}$ are the standard Pochhammer symbols.

 These polynomials obey the three term recurrence relation
\begin{gather*}
 x\mathbf{I}_{n}(x)=\mathbf{I}_{n+1}(x)+(\rho_1-A_n-C_n)\mathbf{I}_{n}(x)+A_{n-1}C_n\mathbf{I}_{n-1}(x)
,\end{gather*}
with recurrence coefficients
\begin{gather}
 A_{n}=\begin{cases}
-\dfrac{\left(n+2 \rho_{2}-2 r_{2}+1\right)\left(n+2 \rho_{2}-2 r_{1}+1\right)}{4(n+g+1)}, & n \text { even}, \vspace{1mm}\\
-\dfrac{(n+1)\left(n-2 r_{1}-2 r_{2}+1\right)}{4(n+g+1)}, & n \text { odd},
\end{cases}\nonumber\\
 C_{n}=\begin{cases}
\dfrac{\left(n+2 \rho_{1}-2 r_{1}+1\right)\left(n+2 \rho_{1}-2 r_{2}+1\right)}{4(n+g+1)}, & n \text { even}, \vspace{1mm}\\
\dfrac{(n+2 g+1)\left(n+2 \rho_{1}+2 \rho_{2}+1\right)}{4(n+g+1)}, & n \text { odd},
\end{cases}\label{C2_RRCBICo}
\end{gather}
where $g=\rho_1+\rho_2-r_1-r_2$.

The general Bannai--Ito polynomials can be presented as the Geronimus transformation of those~$\mathbf{I}_{n}$ with parameter $\rho_1$. The general Bannai--Ito are hence given in terms of general complementary Bannai--Ito as follows:
\begin{gather*}
 \mathbf{B}_{n}(x;\rho_1,\rho_2,r_1,r_2) = \mathbf{I}_{n}(x;\rho_1,\rho_2,r_1,r_2) - A_{n-1}\mathbf{I}_{n-1}(x;\rho_1,\rho_2,r_1,r_2).
\end{gather*}
They obey the three term recurrence relation
\begin{gather}
 x\mathbf{B}_{n}(x)=\mathbf{B}_{n+1}(x)+(\rho_1-A_{n-1}-C_n)\mathbf{B}_{n}(x)+A_{n-1}C_{n-1}\mathbf{B}_{n-1}(x)\label{C2_RRBI},
\end{gather}
where $A_{n}$ and $C_{n}$ are as in \eqref{C2_RRCBICo}.

Dunkl-difference equations are presented for the complementary Bannai--Ito and Bannai--Ito polynomials in \cite{CBI} and \cite{BI}, respectively. It can be seen from formulas \eqref{C2_RRCBICo} that the positivity condition $u_n = A_{n-1}C_n>0$ cannot be achieved for all $n\in\mathbb{N}$ for the general complementary Bannai--Ito. Nevertheless, it is possible to obtain a finite set of $N+1$ orthogonal polynomials by using as truncation conditions $u_0 = A_{-1}C_0=0$ and $u_{N+1}=A_{N}C_{N+1}=0$. The first condition is always respected, but a parametrization dependent on $N$ and its parity is needed for the second condition, and it can be realized in many ways. Suppose that N is even. Then, the truncation conditions $(N+2\rho_2-2r_2+1)=0$, $(N+2\rho_2-2r_1+1)=0$ or $(N+2\rho_1+2\rho_2+2)=0$ all lead to the usual complementary Bannai--Ito polynomials. However, using a parametrization such that \begin{gather}\pp{N+2g+2}=0,\label{C2_trunc_2j}\end{gather} when $N$ is even leads to an admissible truncation and to a different result. One should be careful when choosing a parametrization for this last truncation, since a zero is introduced in the denominator of the coefficients $A_n$ and $C_n$ for $n \sim N / 2$. The parametrization needs to ensure that a cancellation occurs to have finite expressions in those cases. A very similar truncation arises for the general Bannai--Ito but for $N$ odd, and it reads
\begin{gather}(N+2g+1)=0.\label{C2_trunc_2j1}\end{gather}
These truncation conditions lead to a new family of orthogonal polynomials: the para-Bannai--Ito polynomials. Taking $N=2j$ in \eqref{C2_trunc_2j} and $N=2j+1$ in \eqref{C2_trunc_2j1}, both truncation conditions are identical and correspond to
\begin{gather}(j+g+1)=0.\label{C2_trunc_both}\end{gather}
The parametrization then only depends on the parity of $j$.

\section[Para-Bannai--Ito polynomials for N=2j, j even]{Para-Bannai--Ito polynomials for $\boldsymbol{N=2j}$, $\boldsymbol{j}$ even}\label{sec3}

In this section, we obtain the recurrence relation of the para-Bannai--Ito for $N=2j$, $j$ even, by applying the truncation condition \eqref{C2_trunc_both} to the general complementary Bannai--Ito polynomials. An appropriate parametrization is
\begin{gather}
 \rho_1-r_1 = -\frac{j+1}{2}+e_1t,\qquad
 \rho_2-r_2 = -\frac{j+1}{2}+e_2t\label{C3_param_jeven},
\end{gather}
and \eqref{C2_trunc_both} is achieved in the limit $t\to 0$.
\subsection{Recurrence relation}
Inserting \eqref{C3_param_jeven} in \eqref{C2_RRCBICo} and using the change of parameters
\begin{gather}
 \rho_1 = \frac{b-j-1+a}{4},\qquad \rho_2 = \frac{b-j-1-a}{4},\nonumber\\
 \frac{e_1}{e_1+e_2} = \alpha ,\qquad \frac{e_2}{e_1+e_2} = (1-\alpha)\label{C3_change_param},
\end{gather}
it is seen that $A_j^0$ and $C_j^0$ only depend on $e_1$ and $e_2$, through combinations. That can be described in terms of a single parameter $\alpha$ as above. The recurrence relation for the para-Bannai--Ito polynomials $\mathbf{P}_{n}^{(0)}(x;a,b,\alpha,2j)$ reads
\begin{gather}
x\mathbf{P}_{n}^{(0)}(x)=\mathbf{P}_{n+1}^{(0)}(x)+\pp{\frac{b-j-1+a}{4}-A_n^0-C_n^0}\mathbf{P}_{n}^{(0)}(x)+A_{n-1}^0C_n^0\mathbf{P}_{n-1}^{(0)}(x)\label{C3_RRPCBI_2j},\\
 A_{n}^0= \begin{cases}
-\dfrac{1}{4}(n-j-a), & \ n \text { even},\ n\neq j, \\
-\dfrac{1}{4}\dfrac{(n+1)(n-j-b)}{(n-j)}, & n \text { odd},\\
\dfrac{1}{2}(1-\alpha)a, & n \text { even},\ n= j.
\end{cases}\nonumber\\
C_{n}^0= \begin{cases}
\dfrac{1}{4}(n-j+a), & n \text { even},\ n\neq j, \\
\dfrac{1}{4}\dfrac{(n-2j-1)(n-j+b)}{(n-j)}, & n \text { odd},\\
\dfrac{1}{2}\alpha a, & n \text { even},\ n= j.
\end{cases} \label{C3_RRPCBICo_2j}
\end{gather}
It is now clear that the truncation condition $u_0=u_{2j+1}=0$ is achieved. In order to respect the positivity condition $u_n=A_{n-1}^0C_n^0>0$ for $n\in \{1,2,\dots,2j\}$ and obtain a finite set of orthogonal polynomials, one must choose between the two sets of restrictions on the parameters,
\begin{gather}
 a\leq -j-1,\qquad
|b|\leq 1 , \qquad
 0 \leq \alpha \leq 1 , \qquad \text { or } \nonumber\\
 b\geq j, \qquad
|a+1|\leq 1, \qquad
0 \leq \alpha \leq 1.\label{eq:cond_posi}
\end{gather}
It is useful to note that the coefficients of the recurrence relation \eqref{C3_RRPCBI_2j} are persymmetric when~$\alpha=1/2$ (i.e., $u_n = u_{N-n+1}$ and $b_n=b_{N-n}$ where $b_n$ is the coefficient in front of the diagonal term in the recurrence relation). This persymmetry is seen to be tantamount to both conditions $A_n^0+C_n^0=A_{2j-n}^0+C_{2j-n}^0$ and $A_{n-1}^0C_n^0=A_{2j-n}^0C_{2j-n+1}^0$. These conditions are fulfilled when $\alpha=1/2$ because one can verify that $A_n^0=C_{2j-n}^0$ and $C_n^0=A_{2j-n}^0$ which together solve the two conditions. This property will be useful to obtain the orthogonality relation.

\subsection{Dunkl-difference equation}

It is possible to obtain the Dunkl-difference equation for the para-Bannai--Ito polynomials from the one of the complementary Bannai--Ito polynomials via the truncation \eqref{C3_param_jeven}. This procedure is straightforward for this equation since no divergence occur. It reads as
\begin{gather*}
 \mathcal{D}_{\beta} \mathbf{P}_{n}^{(0)}(x)=\Lambda_{n}^{(\beta)} \mathbf{P}_{n}^{(0)}(x),
\end{gather*}
with the following eigenvalues:
\begin{gather*}
 \Lambda_{2n}^{(\beta)} = n(n-j),\qquad \Lambda_{2n+1}^{(\beta)} = n(n+1-j)+\beta.
\end{gather*}
The operator $\mathcal{D}_{\beta}$ is given by
\begin{gather*}
 \mathcal{D}_{\beta} = \mathcal{D}_{0}+\beta \frac{(x-\rho_2)}{2x}(I-R),
\end{gather*}
with
\begin{gather*}
 \mathcal{D}_{0}=A(x) T^{+}+B(x) T^{-}+C(x) R+D(x) T^{+} R-(A(x)+B(x)+C(x)+D(x))I,
\end{gather*}
where $T^{\pm}f(x)=f(x\pm1)$, $Rf(x)=f(-x)$ and
\begin{gather*}
A(x)=\frac{\left(x+\rho_{1}+1\right)\left(x+\rho_{2}+1\right)\left(2 x-2 \rho_{1}-j\right)\left(2 x-2 \rho_{2}-j\right)}{8(x+1)(2 x+1)}, \\
B(x)=\frac{\left(x-\rho_{2}\right)\left(x-\rho_{1}-1\right)\left(2 x+2 \rho_{1}+j\right)\left(2 x+2 \rho_{2}+j\right)}{8 x(2 x-1)}, \\
C(x)=\frac{\left(x-\rho_{2}\right)\left(4 x^{2}+\omega\right)}{8 x}-\frac{\left(x-\rho_{2}\right)\left(x+\rho_{1}+1\right)\left(2 x-2 \rho_{1}-j\right)\left(2 x-2 \rho_{2}-j\right)}{8 x(2 x+1)}-B(x), \\
D(x)=\frac{\rho_{2}\left(x+\rho_{1}+1\right)\left(2 x-2 \rho_{1}-j\right)\left(2 x-2 \rho_{2}-j\right)}{8 x(x+1)(2 x+1)},\\
\omega = (2\rho_1+j)(2\rho_2+j)-4(1+\rho_1)(\rho_1+\rho_2+j).
\end{gather*}
One needs to use the change of parameters \eqref{C3_change_param} in order to obtain the result above. Observe that~$\Lambda_{n}^{(\beta)}=\Lambda_{2j-n}^{(\beta)}$ and that the spectrum is hence degenerate. This means that the para-Bannai--Ito polynomials are bispectral, but not classical.

\subsection{Explicit expression}\label{subsec:Ee}
The explicit expression for the para-Bannai--Ito is obtained by using the limit described after~\eqref{C3_param_jeven} in \eqref{C2_I2n} and \eqref{C2_I2n1}. Let $A_{n,k}^{(I)}$ and $B_{n,k}^{(I)}$ be the coefficients of decomposition of the Wilson polynomials as in \eqref{decompositionwilson} for $\mathbf{I}_{2n}(x)$ and $\mathbf{I}_{2n+1}(x)$. The complementary Bannai--Ito polynomials can then be written as
\begin{gather*}
 \mathbf{I}_{2n}(x) = l_n^{(1)}\sum_{k=0}^{n}A_{n,k}^{(I)}, \qquad \mathbf{I}_{2n+1}(x) = l_n^{(2)}(x-\rho_2)\sum_{k=0}^{n}B_{n,k}^{(I)}.
\end{gather*}
Using this decomposition, the para-Bannai--Ito polynomials can be obtained by taking the limit directly in those coefficients,
\begin{gather*}
 \mathbf{P}_{2n}^{(0)}(x) = \kappa_n^{(1)}\sum_{k=0}^{n} \lim_{t\to 0} A_{n,k}^{(I)}= \kappa_n^{(1)}\sum_{k=0}^{n} A_{n,k},\\ \mathbf{P}_{2n+1}^{(0)}(x) = \kappa_n^{(2)}(x-\rho_2)\sum_{k=0}^{n}\lim_{t\to 0} B_{n,k}^{(I)}= \kappa_n^{(2)}\pp{x-\frac{b-j-1-a}{4}}\sum_{k=0}^{n} B_{n,k},
\end{gather*}
with the right normalization to make them monic. The summands and the renormalization constants are found to be
\begin{gather*}% \label{C_Ank}\label{C_kap1}
A_{n,k}=\begin{cases}
\dfrac{(-n)_{k}(n-j)_{k}\bigl(\frac{b-j-1-a}{4}+x\bigr)_{k}\bigl(\frac{b-j-1-a}{4}-x\bigr)_{k}}{(1_{k}\bigl(\frac{1+b-j}{2}\bigr)_{k}
\bigl(-\frac{j+a}{2}\bigr)_{k}\bigl(-\frac{j}{2}\bigr)_{k}},\\
\hspace*{100mm} k\leq \frac{j}{2},\ k\leq n, \\[0.3cm]
\pp{\dfrac{1}{1-\alpha}}\dfrac{(-n)_{k}\bigl(\frac{b-j-1-a}{4}+x\bigr)_{k}\bigl(\frac{b-j-1-a}{4}-x\bigr)_{k}}{(1)_{k}\bigl(\frac{1+b-j}{2}\bigr)_{k}
\bigl(-\frac{j+a}{2}\bigr)_{k}}\dfrac{(n-j)_{j-n}(1)_{n+k-j-1}}{\bigl(-\frac{j}{2}\bigr)_{\frac{j}{2}}(1)_{k-j/2-1}}, \\
\hspace*{100mm} k> \frac{j}{2},\ k\leq n,
\end{cases}\\
\kappa_{n}^{(1)}= \begin{cases}
\dfrac{\bigl(\frac{1+b-j}{2}\bigr)_{n}\bigl(-\frac{j+a}{2}\bigr)_{n}\bigl(-\frac{j}{2}\bigr)_{n}}{(n-j)_{n}}, & n\leq \frac{j}{2}, \\[0.3cm]
(1-\alpha)\dfrac{\bigl(\frac{1+b-j}{2}\bigr)_{n}\bigl(-\frac{j+a}{2}\bigr)_{n}\bigl(-\frac{j}{2}\bigr)_{\frac{j}{2}}(1)_{n-j/2-1}}{(n-j)_{j-n}(1)_{2n-j-1}}, & n> \frac{j}{2},
\end{cases}
\end{gather*}
and
\begin{gather*}%\label{C_Bnk}\label{C_kap2}
B_{n,k}= \begin{cases}
\dfrac{(-n)_{k}(n+1-j)_{k}\bigl(\frac{b-j+3-a}{4}+x\bigr)_{k}\bigl(\frac{b-j+3-a}{4}-x\bigr)_{k}}{(1)_{k}\bigl(\frac{3+b-j}{2}\bigr)_{k}\bigl(\frac{2-j-a}{2}\bigr)_{k}
\bigl(\frac{2-j}{2}\bigr)_{k}}, & k< \frac{j}{2}, \vspace{1mm}\\
\!\left(\!\dfrac{1}{1-\alpha}\!\right)\!\dfrac{(-n)_{k}\bigl(\frac{b-j+3-a}{4}+x\bigr)_{k}
\bigl(\frac{b-j+3-a}{4}-x\bigr)_{k}}{(1)_{k}\bigl(\frac{3+b-j}{2}\bigr)_{k}\bigl(\frac{2-j-a}{2}\bigr)_{k}}\dfrac{(n+1-j)_{j-n-1}(1)_{n+k-j}}
{\bigl(\frac{2-j}{2}\bigr)_{\frac{j-2}{2}}(1)_{k-j/2}},\!\! & k\geq \frac{j}{2},
\end{cases}\\
\kappa_{n}^{(2)}= \begin{cases}
\dfrac{\bigl(\frac{3+b-j}{2}\bigr)_{n}\bigl(\frac{2-j-a}{2}\bigr)_{n}\bigl(\frac{2-j}{2}\bigr)_{n}}{(n+1-j)_{n}}, & n< \frac{j}{2}, \\[0.3cm]
(1-\alpha)\dfrac{\bigl(\frac{3+b-j}{2}\bigr)_{n}\bigl(\frac{2-j-a}{2}\bigr)_{n}\bigl(\frac{2-j}{2}\bigr)_{\frac{j-2}{2}}(1)_{n-j/2}}{(n+1-j)_{j-n-1}(1)_{2n-j}}, & n\geq \frac{j}{2}.
\end{cases}
\end{gather*}
In terms of terminating hypergeometric series, for even degrees, we have
\begin{gather*}
 \frac{1}{\kappa_n^{(1)}}\mathbf{P}_{2n}^{(0)}(x) = \pFq{4}{3}{-n,n-j,\frac{b-j-1-a}{4}+x,\frac{b-j-1-a}{4}-x}{\frac{1+b-j}{2},-\frac{j+a}{2},-\frac{j}{2}}{1} ,
\end{gather*}
if $n\leq \frac{j}{2}$, and
\begin{align*}
 \frac{1}{\kappa_n^{(1)}}\mathbf{P}_{2n}^{(0)}(x) ={}&\pFq{4}{3}{-n,n-j,\frac{b-j-1-a}{4}+x,\frac{b-j-1-a}{4}-x}{\frac{1+b-j}{2},-\frac{j+a}{2},-\frac{j}{2}}{1}\nonumber\\ &+\frac{(n-j)_{j-n}(-n)_{j/2+1}\bigl(\frac{b-j-1-a}{4}+x\bigr)_{j/2+1}\bigl(\frac{b-j-1-a}{4}-x\bigr)_{j/2+1}(1)_{n-j/2}} {(1-\alpha)\bigl(-\frac{j}{2}\bigr)_{j/2}(1)_{j/2+1}\bigl(\frac{1+b-j}{2}\bigr)_{j/2+1}\bigl(-\frac{j+a}{2}\bigr)_{j/2+1}}\nonumber\\
 &\phantom{+}{}\times \pFq{4}{3}{\frac{j+2}{2}-n,n-\frac{j-2}{2},\frac{b+j+3-a}{4}+x,\frac{b+j+3-a}{4}-x}{\frac{3+b}{2},\frac{2-a}{2},\frac{4+j}{2}}{1}
\end{align*}
if $n> \frac{j}{2}$.
For odd degrees, we have
\begin{gather*}
 \frac{1}{\kappa_n^{(2)}\bigl(x-\frac{b-j-1-a}{4}\bigr)}\mathbf{P}_{2n+1}^{(0)}(x) = \pFq{4}{3}{-n,n+1-j,\frac{b-j+3-a}{4}+x,\frac{b-j+3-a}{4}-x}{\frac{3+b-j}{2},\frac{2-j-a}{2},\frac{2-j}{2}}{1},
\end{gather*}
if $n< \frac{j}{2}$, and
\begin{gather*}
 \frac{1}{\kappa_n^{(2)}\bigl(x-\frac{b-j-1-a}{4}\bigr)}\mathbf{P}_{2n+1}^{(0)}(x) =\pFq{4}{3}{-n,n+1-j,\frac{b-j+3-a}{4}+x,\frac{b-j+3-a}{4}-x}{\frac{3+b-j}{2},\frac{2-j-a}{2},\frac{2-j}{2}}{1}\nonumber\\[0.3cm] \qquad+\frac{(n+1-j)_{j-n-1}(-n)_{j/2}\bigl(\frac{b-j+3-a}{4}+x\bigr)_{j/2}\bigl(\frac{b-j+3-a}{4}-x\bigr)_{j/2}(1)_{n-j/2}} {(1-\alpha)\bigl(\frac{2-j}{2}\bigr)_{j/2-1}(1)_{j/2}\bigl(\frac{3+b-j}{2}\bigr)_{j/2}\bigl(\frac{2-j-a}{2}\bigr)_{j/2}}\nonumber\\[0.3cm]
 \qquad \qquad\times\pFq{4}{3}{\frac{j}{2}-n,n-\frac{j-2}{2},\frac{b+j+3-a}{4}+x,\frac{b+j+3-a}{4}-x}{\frac{3+b}{2},\frac{2-a}{2},\frac{4+j}{2}}{1}
\end{gather*}
if $n\geq \frac{j}{2}$.

\subsection{Orthogonality relation}\label{C3_OR}
To obtain the orthogonality relation, we first need to determine the grid or orthogonality lattice. It is provided by the eigenvalues of the Jacobi matrix or, equivalently, by the zeroes of the characteristic polynomial \smash{$\mathbf{P}_{2j+1}^{(0)}(x)$}. We can define \smash{$\mathbf{P}_{2j+1}^{(0)}(x)$} using the recurrence relation \eqref{C3_RRPCBI_2j} and inserting the expression of \smash{$\mathbf{P}_{2j}^{(0)}(x)$} and \smash{$\mathbf{P}_{2j-1}^{(0)}(x)$} given in Section~\ref{subsec:Ee}. An extraction of the zeroes can be done with the help of the Saalsch\"utz summation formula to obtain
\begin{gather*}
\mathbf{P}_{2j+1}^{(0)}(x) = \prod_{k=0}^{j}\pp{x+(-1)^k\pp{\frac{2k-j-b+a}{4}}+\frac{1}{4}}\\
\hphantom{\mathbf{P}_{2j+1}^{(0)}(x) =}{}
\times \prod_{k=0}^{j-1}\pp{x+(-1)^k\pp{\frac{2k-j-b-a}{4}}+\frac{1}{4}}.
\end{gather*}
The set of $2j+1$ para-Bannai--Ito polynomials will thus be orthogonal on a Bannai--Ito bi-lattice with $2j+1$ grid points.
\begin{gather*}
 x_{2s} = -(-1)^s\pp{\frac{2s-j-b+a}{4}}-\frac{1}{4},\qquad s\in \{0,1,\dots,j\},\\
 x_{2s+1} = -(-1)^s\pp{\frac{2s-j-b-a}{4}}-\frac{1}{4},\qquad s\in \{0,1,\dots,j-1\}.
\end{gather*}
A depiction of the grid is given in Appendix~\ref{secappC}. This Bannai--Ito bi-lattice can also be seen as a linear quadri-lattice
\begin{gather*}
 x_{4s} = -\pp{\frac{4s-j-b+a}{4}}-\frac{1}{4},\qquad s\in \{0,1,\dots,j/2\},\\
 x_{4s+1} = -\pp{\frac{4s-j-b-a}{4}}-\frac{1}{4},\qquad s\in \{0,1,\dots,j/2-1\},\\
 x_{4s+2} = \pp{\frac{4s+2-j-b+a}{4}}-\frac{1}{4},\qquad s\in \{0,1,\dots,j/2-1\},\\
 x_{4s+3} = \pp{\frac{4s+2-j-b-a}{4}}-\frac{1}{4},\qquad s\in \{0,1,\dots,j/2-1\}.
\end{gather*}
From the theory of orthogonal polynomials \cite{CHIbook}, the weights are given by
\begin{gather}
 w_{s}=\frac{u_{1} \cdots u_{2j}}{\mathbf{P}_{2j}^{(0)}(x_s) \mathbf{P}_{2j+1}^{(0)\prime}(x_s)}, \qquad s=0,1, \ldots, 2j\label{C3_full_weight},
\end{gather}
in terms of which the orthogonality relation of the para-Bannai--Ito polynomials will read
\begin{gather}
 \sum_{s=0}^{2j}w_{s} \mathbf{P}_{n}^{(0)}(x_s)\mathbf{P}_{m}^{(0)}(x_s) = u_{1} \cdots u_{n} \delta_{nm}\label{C3_ortho}.
\end{gather}
In \cite{persymmetric}, a method to derive the weight for persymmetric orthogonal polynomials was presented together with a procedure to also achieve that after an isospectral deformation. This is exactly the situation that the para-Bannai--Ito polynomials present. Now, the usual expression
\begin{gather}
 \mathbf{P}_{N}^{(0)}(x_s) = (-1)^{N+s}\sqrt{u_{1} \cdots u_{N}}\label{C3_simp_p2j_xs}
\end{gather}
for $\mathbf{P}_{N}^{(0)}(x_s)$ is predicated on the eigenvalues of the Jacobi matrix ordered in an increasing fashion which is not the case in the presentation of the bi-Bannai--Ito grid given above. Nonetheless, it is possible to alter the increasing order condition and to require instead that between any two zeroes of even index, there must be a zero of odd index and vice versa. Using this modified condition, the interlacing zeroes property and the persymmetry, when $N=2j$, we still find for~\smash{$\mathbf{P}_{N}^{(0)}(x_s)$} the expression \eqref{C3_simp_p2j_xs}. In the case of the para-Bannai--Ito polynomials, positivity ensures that the modified condition for the eigenvalues is achieved. Recalling that the polynomials are persymmetric, we then have that for $\alpha=1/2$,
\begin{gather}
 \tilde{w}_{s} = \frac{(-1)^{s}\sqrt{u_{1} \cdots u_{N}}}{\mathbf{P}_{2j+1}^{(0)\prime}(x_s)} \label{C3_redu_weight},
\end{gather}
with the $u_n$ also evaluated at $\alpha=1/2$. A direct computation gives
\begin{gather*}
 \tilde{w}_{4s} = \frac{h_{2j}}{(-s)_{s}(1)_{j/2-s}\bigl(-s-\frac{a}{2}\bigr)_{j/2}\bigl(s+\frac{1+a-b-j}{2}\bigr)_{j/2}\bigl(s+\frac{1-b-j}{2}\bigr)_{j/2}},\\
 \tilde{w}_{4s+1} = -\frac{h_{2j}}{(-s)_{s}(1)_{j/2-s-1}\bigl(-s+\frac{a}{2}\bigr)_{j/2+1}\bigl(s+\frac{1-a-b-j}{2}\bigr)_{j/2}\bigl(s+\frac{1-b-j}{2}\bigr)_{j/2}},\\
 \tilde{w}_{4s+2} = -\frac{h_{2j}}{(-s)_{s}(1)_{j/2-s-1}\bigl(-s-\frac{a}{2}\bigr)_{j/2}\bigl(s+\frac{1+a-b-j}{2}\bigr)_{j/2+1}\bigl(s+\frac{1-b-j}{2}\bigr)_{j/2}},\\
 \tilde{w}_{4s+3} = \frac{h_{2j}}{(-s)_{s}(1)_{j/2-s-1}\bigl(-s+\frac{a}{2}\bigr)_{j/2}\bigl(s+\frac{1-a-b-j}{2}\bigr)_{j/2}\bigl(s+\frac{1-b-j}{2}\bigr)_{j/2+1}},
\end{gather*}
with
\begin{gather*}
 h_{2j} = \sqrt{u_{1} \cdots u_{2j}}= \frac{(1)_{j}\bigl(\frac{b+1-j}{2}\bigr)_{j}\bigl(-\frac{j+a}{2}\bigr)_{j}}{2^{2j}\bigl(\frac{1-j}{2}\bigr)_{j}}.
\end{gather*}
It was also shown in \cite{persymmetric} that after an isospectral deformation, the new weights are related to the persymmetric weights by a multiplicative factor dependent on $\alpha$ in our case,
\begin{gather*}
 w_s = C\pp{1+\beta(-1)^s}\tilde{w}_s.
\end{gather*}
Using \eqref{C3_full_weight} and \eqref{C3_redu_weight} for the first few $j$, we can solve for $C$ and $\beta$. One can verify that $C=1$ and $\beta = 2\alpha-1$ which gives
\begin{gather*}
w_{2 s} =2\alpha \tilde{w}_{2 s},\qquad
w_{2 s+1} =2 (1-\alpha) \tilde{w}_{2 s+1}.
\end{gather*}
We also have
\begin{gather*}
 \sum_{s=0}^{j} w_{2 s}=\alpha, \qquad \sum_{s=0}^{j} w_{2 s+1}=1-\alpha,
\end{gather*}
which generalize a known result for mirror-symmetric Jacobi matrices. Finally, we give a general expression for $h_n = \sqrt{u_{1} \cdots u_{n}}$ which can be squared and inserted in \eqref{C3_ortho} to complete the characterization of the para-Bannai--Ito for $N=2j$ with $j$ even,
\begin{gather*}
 h_{2n} = \frac{\sqrt{(-j)_{n}(1)_{n}\bigl(\frac{1+b-j}{2}\bigr)_{n}\bigl(\frac{1-b-j}{2}\bigr)_{n}\bigl(-\frac{j+a}{2}\bigr)_{n}
 \bigl(\frac{2+a-j}{2}\bigr)_{n}}}{2^{2n}\bigl(\frac{1-j}{2}\bigr)_{n}}\times \begin{cases}
1, & n<\frac{j}{2}, \\[0.1cm]
\sqrt{2\alpha}, & n=\frac{j}{2},\\[0.1cm]
2\sqrt{\alpha(1-\alpha)}, & n>\frac{j}{2},
\end{cases} \\
h_{2n+1} = \frac{\sqrt{(-1)(-j)_{n+1}(1)_{n}\bigl(\frac{1+b-j}{2}\bigr)_{n+1}\bigl(\frac{1-b-j}{2}\bigr)_{n} \bigl(-\frac{j+a}{2}\bigr)_{n+1}\bigl(\frac{2+a-j}{2}\bigr)_{n}}}{2^{2n+1}\sqrt{\bigl(\frac{1-j}{2}\bigr)_{n}\bigl(\frac{1-j}{2}\bigr)_{n+1}}}\\ \phantom{h_{2n+1} = }{}\times \begin{cases}
1, & n<\frac{j}{2}, \\[0.1cm]
2\sqrt{\alpha(1-\alpha)}, & n\geq \frac{j}{2}.
\end{cases}\nonumber
\end{gather*}

\section[Para-Bannai--Ito for N=2j+1, j even]{Para-Bannai--Ito for $\boldsymbol{N=2j+1}$, $\boldsymbol{j}$ even}\label{sec4}
The general complementary Bannai--Ito and the general Bannai--Ito polynomials are related by a Geronimus transformation with parameter $\rho_1$. The Geronimus transformation and the truncation commute, and therefore, performing the transformations on the para-Bannai--Ito polynomials with $N=2j+1$, $j$ even, we obtain a set of $N+1=2j+2$ polynomials that correspond to the truncation \eqref{C3_param_jeven} of the general Bannai--Ito polynomials. The eigenvalue ${x_s = \rho_1}$ is then added to the spectrum. The treatment of the case $N=2j+1$, $j$ even, is analogous to the case $N=2j$, $j$ even. The results are presented below.

\subsection{Recurrence relation}
Using the truncation \eqref{C3_param_jeven} in \eqref{C2_RRBI} gives the recurrence relation for $\mathbf{P}_{n}^{(1)}\pp{x;a,b,\alpha,2j+1}$
\begin{gather}%\label{C4_RRPCBICo_2j}
 x\mathbf{P}_{n}^{(1)}(x)=\mathbf{P}_{n+1}^{(1)}(x)+\pp{\frac{b-j-1+a}{4}-A_n^1-C_n^1}\mathbf{P}_{n}^{(1)}(x)+A_{n-1}^1C_n^1\mathbf{P}_{n-1}^{(1)}(x)\label{C4_RRPBI_2j1},\\
 A_{n}^1= \begin{cases}
\dfrac{1}{4}(n-j+a), & n \text { even},\ n\neq j, \\[0.3cm]
\dfrac{1}{4}\dfrac{(n-2j-1)(n-j+b)}{(n-j)}, & n \text{ odd},\\[0.3cm]
\dfrac{1}{2}\alpha a, & n \text { even},\ n= j,
\end{cases}\nonumber\\
C_{n}^1= \begin{cases}
-\dfrac{1}{4}\dfrac{n(n-j-1-b)}{(n-j-1)}, & n \text{ even}, \\[0.3cm]
-\dfrac{1}{4}(n-j-1-a), & n \text { odd},\ n\neq j+1, \\[0.3cm]
\dfrac{1}{2}(1-\alpha) a, & n \text { odd},\ n= j+1.
\end{cases}\nonumber
\end{gather}
The positivity condition is achieved if one chooses the parameters such that
\begin{gather*}
 |a|\geq j+1,\qquad
 |b|\leq 1,\qquad
 0 \leq \alpha \leq 1 .
\end{gather*}
Again, the polynomials are persymmetric if $\alpha=1/2$.
\subsection{Dunkl-difference equation}
Since no divergences appear in the Dunkl-difference equation of the general Bannai--Ito polynomials under the para truncation, it follows that the para-Bannai--Ito polynomials are solutions to the same equation except for a change of variable. It reads as
\begin{gather*}
L \mathbf{P}_{n}^{(1)}(x)=\lambda_{n} \mathbf{P}_{n}^{(1)}(x),
\end{gather*}
with the eigenvalues
\begin{gather*}
 \lambda_{2n} = n, \qquad \lambda_{2n+1} = \pp{j-n}.
\end{gather*}
The Dunkl-difference operator is
\begin{gather*}
 L=F(x)(I-R)+G(x)\bigl(T^{+} R-I \bigr),
\end{gather*}
where
\begin{gather*}
 G(x) = \frac{(4x+1+a-b-j)(4x+1-a-b-j)}{16(2x+1)},\\ F(x) = \frac{(4x+1+j+a-b)(4x+1+j-a-b)}{32x}.
\end{gather*}
\subsection{Explicit expression}
Since the para-Bannai--Ito for $N=2j$ and $N=2j+1$ are related by a Geronimus transformation, we obtain directly the expression for the case $N=2j+1$ with $j$ even as a combination of two para-Bannai--Ito polynomials with $N=2j$. The explicit expression is
\begin{gather*}
 \mathbf{P}_{n}^{(1)}(x) = \mathbf{P}_{n}^{(0)}(x)-C_n^{1}\mathbf{P}_{n-1}^{(0)}(x),
\end{gather*}
where $C_n^{1}$ is as in \eqref{C4_RRPBI_2j1}.
\subsection{Orthogonality relation}
We know that the Geronimus transformation will only add $x_s = \rho_1$ to the spectrum provided by the zeros of \smash{$\mathbf{P}_{2j+1}^{(0)}$}. It turns out that this fits with the grid of Section \ref{C3_OR} and that its points are then given by
\begin{gather}
 x_{2s} = -(-1)^s\pp{\frac{2s-j-b+a}{4}}-\frac{1}{4},\qquad s\in \{0,1,\dots,j\}\label{C4_spec_1},\\
 x_{2s+1} = -(-1)^s\pp{\frac{2s-j-b-a}{4}}-\frac{1}{4},\qquad s\in \{0,1,\dots,j\}\label{C4_spec_2},
\end{gather}
or again as a linear quadri-lattice
\begin{gather*}
 x_{4s} = -\pp{\frac{4s-j-b+a}{4}}-\frac{1}{4},\qquad s\in \{0,1,\dots,j/2\},\\
 x_{4s+1} = -\pp{\frac{4s-j-b-a}{4}}-\frac{1}{4},\qquad s\in \{0,1,\dots,j/2\},\\
 x_{4s+2} = \pp{\frac{4s+2-j-b+a}{4}}-\frac{1}{4},\qquad s\in \{0,1,\dots,j/2-1\},\\
 x_{4s+3} = \pp{\frac{4s+2-j-b-a}{4}}-\frac{1}{4},\qquad s\in \{0,1,\dots,j/2-1\}.
\end{gather*}
The orthogonality relation is again of the form
\begin{gather*}
 \sum_{s=0}^{2j+1}w_{s} \mathbf{P}_{n}^{(1)}(x_s)\mathbf{P}_{m}^{(1)}(x_s) = u_{1} \cdots u_{n} \delta_{nm}\label{C4_ortho},
\end{gather*}
where
\begin{gather*}
 w_{s}=\frac{u_{1} \cdots u_{2j+1}}{\mathbf{P}_{2j+1}^{(1)}(x_s) \mathbf{P}_{2j+2}^{(1)\prime}(x_s)}, \qquad s=0,1, \ldots, 2j+1\label{C4_full_weight}.
\end{gather*}
The computation in the persymmetric case followed by the generalization with an isospectral deformation gives the weights
\begin{gather*}
 w_{4s} = \frac{2 \alpha h_{2j+1}}{(-s)_{s}(1)_{j/2-s}\bigl(-s-\frac{a}{2}\bigr)_{j/2+1}\bigl(s+\frac{1+a-b-j}{2}\bigr)_{j/2}\bigl(s+\frac{1-b-j}{2}\bigr)_{j/2}},\\
 w_{4s+1} = -\frac{2 (1-\alpha) h_{2j+1}}{(-s)_{s}(1)_{j/2-s}\bigl(-s+\frac{a}{2}\bigr)_{j/2+1}\bigl(s+\frac{1-a-b-j}{2}\bigr)_{j/2}\bigl(s+\frac{1-b-j}{2}\bigr)_{j/2}},\\
 w_{4s+2} = -\frac{2 \alpha h_{2j+1}}{(-s)_{s}(1)_{j/2-s-1}\bigl(-s-\frac{a}{2}\bigr)_{j/2}\bigl(s+\frac{1+a-b-j}{2}\bigr)_{j/2+1}\bigl(s+\frac{1-b-j}{2}\bigr)_{j/2+1}},\\
 w_{4s+3} = \frac{2 (1-\alpha)h_{2j+1}}{(-s)_{s}(1)_{j/2-s-1}\bigl(-s+\frac{a}{2}\bigr)_{j/2}\bigl(s+\frac{1-a-b-j}{2}\bigr)_{j/2+1}\bigl(s+\frac{1-b-j}{2}\bigr)_{j/2+1}},
\end{gather*}
with
\begin{gather*}
 h_{2j+1} = \frac{(1)_{j}\bigl(\frac{1-b-j}{2}\bigr)_{j}\bigl(-\frac{j+a}{2}\bigr)_{j+1}}{2^{2j+1}\bigl(\frac{1-j}{2}\bigr)_{j}}.
\end{gather*}
To complete the characterization for $N=2j+1$ with $j$ even, we have that the general expressions for the $h_n = \sqrt{u_{1} \cdots u_{n}}$ with $u_n = A_{n-1}^1C_{n}^1$ are
\begin{gather*}
 h_{2n} = \frac{\sqrt{(-j)_{n}(1)_{n}\bigl(\frac{1-b-j}{2}\bigr)_{n}\bigl(\frac{1+b-j}{2}\bigr)_{n} \bigl(-\frac{j+a}{2}\bigr)_{n}\bigl(-\frac{j-a}{2}\bigr)_{n}}}{2^{2n}\bigl(\frac{1-j}{2}\bigr)_{n}}\times \begin{cases}
1, & n\leq \frac{j}{2}, \\[0.1cm]
2\sqrt{\alpha(1-\alpha)}, & n> \frac{j}{2},
\end{cases}\\
h_{2n+1} = \frac{\sqrt{(-1)(-j)_{n}(1)_{n}\bigl(\frac{1-b-j}{2}\bigr)_{n}\bigl(\frac{1+b-j}{2}\bigr)_{n} \bigl(-\frac{j+a}{2}\bigr)_{n+1}\bigl(-\frac{j-a}{2}\bigr)_{n+1}}}{2^{2n+1}\bigl(\frac{1-j}{2}\bigr)_{n}}\\ \phantom{h_{2n+1} =}{}\times \begin{cases}
1, & n<\frac{j}{2}, \\[0.1cm]
2\sqrt{\alpha(1-\alpha)}, & n\geq \frac{j}{2}.
\end{cases}
\end{gather*}

\section[Para-Bannai--Ito polynomials as a q to -1 limit of the q-para-Racah polynomials]{Para-Bannai--Ito polynomials as a $\boldsymbol{q\to-1}$ limit\\ of the $\boldsymbol{q}$-para-Racah polynomials}\label{sec5}

In \cite{qpararacah}, the $q$-para-Racah polynomials were introduced as a $q$ generalization of the para-Racah polynomials. In this section, we show that the para-Bannai--Ito polynomials can be obtained as a $q\to -1$ limit of the $q$-para-Racah polynomials, reinforcing their position in the family of para polynomials. The monic $q$-para-Racah polynomials for $N=2j+1$ satisfy the recurrence relation
\begin{gather*}
 x\mathbf{R}_{n}(x)=\mathbf{R}_{n+1}(x)+\pp{\frac{c+c^{-1}}{2}-A_n^R-C_n^R}\mathbf{R}_{n}(x)+A_{n-1}^RC_n^R\mathbf{R}_{n-1}(x)\label{Cr_RRqPR_2j1},
\end{gather*}
with
\begin{gather*}
A_{n}^R= \begin{cases}
\dfrac{\left(1-c d q^{n}\right)\left(d-c q^{n-j}\right)\left(1-q^{n-2 j-1}\right)}{2c d\left(1-q^{2 n-2 j-1}\right)\left(1+q^{n-j}\right)}, & n \neq j, \\[0.2cm]
\dfrac{\alpha\left(1-c d q^{j}\right)(d-c)\left(1-q^{-j-1}\right)}{2c d\left(1-q^{-1}\right)}, & n=j,
\end{cases}\\
C_{n}^R= \begin{cases}
\dfrac{\left(1-q^{n}\right)\left(c-d q^{n-j-1}\right)\left(c d-q^{n-2 j-1}\right)}{2c d\left(1+q^{n-j-1}\right)\left(1-q^{2 n-2 j-1}\right)}, & n \neq j+1, \\[0.2cm]
\dfrac{(1-\alpha)\left(1-q^{j+1}\right)(c-d)\left(c d-q^{-j}\right)}{2c d(1-q)}, & n=j+1.
\end{cases}
\end{gather*}
In order to recover the para-Bannai--Ito polynomials, one needs to use the parametrization
\begin{gather}\label{C5_para_qPR_to_PBI}
 q=-{\rm e}^\varepsilon, \qquad
 c={\rm i}\expp{\varepsilon\bb{\frac{a+b-j}{2}}}, \qquad
 d={\rm i}\expp{\varepsilon\bb{\frac{-a+b-j}{2}}}.
\end{gather}
After redefining the resulting family of polynomials through the following affine transformation of the variable
\begin{gather*}
 x\to 2{\rm i}\varepsilon\pp{x+\frac{1}{4}},
\end{gather*}
we have
\begin{gather}
 x\mathbf{\tilde{R}}_{n}(x)=\mathbf{\tilde{R}}_{n+1}(x)+\pp{\frac{c+c^{-1}}{4{\rm i}\varepsilon}-\frac{1}{4}-\frac{A_n^R}{2{\rm i}\varepsilon}-\frac{C_n^R}{2{\rm i}\varepsilon}}\mathbf{\tilde{R}}_{n}(x)+\frac{A_{n-1}^R}{2{\rm i}\varepsilon}\frac{C_n^R}{2{\rm i}\varepsilon}\mathbf{\tilde{R}}_{n-1}(x)\label{New_C5_RRqPR_2j1}.
\end{gather}
Using the parametrization \eqref{C5_para_qPR_to_PBI} and taking the limit $\varepsilon\to0$, equivalent to $q\to -1$, one sees that~\eqref{New_C5_RRqPR_2j1} goes into the recurrence relation of the para-Bannai--Ito polynomials with $N=2j+1$ and $j$ even. While this was not discussed in \cite{qpararacah}, the $N=2j+1$ and $N=2j$ cases of the $q$-para-Racah polynomials are related by Christoffel and Geronimus transformations, which is the analog of the relation presented before. So the limit process also works for the $q$-para-Racah with $N=2j$, and gives the para-Bannai--Ito with $N=2j$ and $j$ even and the recurrence relation~\eqref{C3_RRPCBICo_2j}.

\section{Special cases}
\subsection[Reduction to a simple lattice and connection with the dual -1 Hahn polynomials]{Reduction to a simple lattice\\ and connection with the dual $\boldsymbol{-1}$ Hahn polynomials}

Consider the case $\alpha=1/2$. For $b=0$, it is observed that the spectrum of the para-Bannai--Ito (i.e., \eqref{C4_spec_1} and \eqref{C4_spec_2}) reduces to a single Bannai--Ito lattice (two linear lattices separated by a gap). This reduction can be observed in the appendix. This restricted grid can be expressed by the sequence
\begin{gather*}
 2 x_k = \frac{(-1)^k}{4}\pp{\frac{k}{2}+\frac{1}{4}-\frac{a+j+1}{2}}-\frac{1}{4},\qquad k\in \{0,1,\dots,N\}.
\end{gather*}
In this setting, we can connect this special case directly to the general complementary Bannai--Ito and to the general Bannai--Ito polynomials
\begin{gather*}
 \mathbf{P}_{n}^{(0)}\pp{x;a,0,\frac{1}{2},2j} = \mathbf{I}_{n}\pp{x;-\frac{j+1-a}{4},-\frac{j+1+a}{4},\frac{j+1-a}{4},\frac{j+1+a}{4}},\\
 \mathbf{P}_{n}^{(1)}\pp{x;a,0,\frac{1}{2},2j+1} = \mathbf{B}_{n}\pp{x;-\frac{j+1-a}{4},-\frac{j+1+a}{4},\frac{j+1-a}{4},\frac{j+1+a}{4}},
\end{gather*}
where $\mathbf{I}_{n}$ and $\mathbf{B}_{n}$ are as in Section \ref{C2_GCBI_GBI}.

 This special case is also connected with the dual $-1$ Hahn polynomials $\mathbf{\hat{R}}_{n}(x;\alpha,\beta,N)$ \cite{D-1H} via the relations
\begin{gather*}
 \mathbf{P}_{n}^{(0)}\pp{x;a,0,\frac{1}{2},2j} = \frac{1}{2^{3n}}\mathbf{\hat{R}}_{n}\bigl(-\bigl[2^3x+1\bigr];j-a,j-a,2j\bigr),\\
 \mathbf{P}_{n}^{(1)}\pp{x;a,0,\frac{1}{2},2j+1} = \frac{1}{2^{3n}}\mathbf{\hat{R}}_{n}\bigl(\bigl[2^3x+1\bigr];a-j-1,a-j-1,2j+1\bigr).
\end{gather*}

\subsection{Reduction to a simple linear lattice}
The simple Bannai--Ito grid can be further specialized by taking $a=-j-1$. This takes it into a simple linear lattice
\begin{gather*}
 x_k = j-k,\qquad k\in \{0,1,\dots,N\},
\end{gather*}
after a relabelling of the grid points. The polynomials are then related to shifted monic Krawtchouk polynomials $\mathbf{K}_{n}(x;p,N)$ \cite{koek}
\begin{gather*}
 \mathbf{P}_{n}^{(0)}\pp{x;-j-1,0,\frac{1}{2},2j} = \frac{1}{2^n}\mathbf{K}_{n}\pp{2x+j;\frac{1}{2},2j},\\
 \mathbf{P}_{n}^{(1)}\pp{x;-j-1,0,\frac{1}{2},2j+1} = \frac{1}{2^n}\mathbf{K}_{n}\pp{2x+j+1;\frac{1}{2},2j+1}.
\end{gather*}

\section{Conclusion}
To summarize, the para-Bannai--Ito polynomials have been introduced and characterized. They are obtained by a $q\to -1$ limit of the $q$-para-Racah polynomials, but also as a special truncation of the general Bannai--Ito and complementary Bannai--Ito polynomials. Their explicit expressions in terms of hypergeometric series was derived as well as their recurrence relation and Dunkl-difference equation. It was shown that the para-Bannai--Ito polynomials are orthogonal on a finite bi-Bannai--Ito grid (linear quadri-lattice). The para-Bannai--Ito have a connection to the dual $-1$ Hahn polynomials (orthogonal on single Bannai--Ito grid or linear bi-lattice) and to the shifted Krawtchouk polynomials (orthogonal on a linear lattice).

The Jacobi matrix of the para-Bannai--Ito is an isospectral deformation of the persymmetric one when $\alpha=1/2$. Such matrices have been used in the design of spin chains with fractional revival. It would be of interest to examine if models based on the para-Bannai--Ito polynomials exhibit similar properties.

\appendix
\section[Para-Bannai--Ito for N=2j, j odd]{Para-Bannai--Ito for $\boldsymbol{N=2j}$, $\boldsymbol{j}$ odd}\label{sec:appA}
\textit{Truncation from the general complementary Bannai--Ito polynomials}
\begin{gather*}
 \rho_1+\rho_2 = -\frac{j+1}{2}+e_1t,\qquad
 -r_1-r_2 = -\frac{j+1}{2}+e_2t\label{CA_param_jodd}.
\end{gather*}
The para-Bannai--Ito polynomials are obtained through a $t\to0$ limit.
\begin{gather}
 \rho_1 = \frac{b-j-1+a}{4},\qquad -r_1 = \frac{b-j-1-a}{4},\nonumber\\ \frac{e_1}{e_1+e_2} = \alpha, \qquad \frac{e_2}{e_1+e_2} = (1-\alpha)\label{CA_change_param}.
\end{gather}
\textit{Recurrence relation}
\begin{gather*}%\label{CA_RRPBICo_2j}\label{CA_RRPBI_2j}
 x\mathbf{P}_{n}^{(2)}(x)=\mathbf{P}_{n+1}^{(2)}(x)+\pp{\frac{b-j-1+a}{4}-A_n^2-C_n^2}\mathbf{P}_{n}^{(2)}(x)+A_{n-1}^2C_n^2\mathbf{P}_{n-1}^{(2)}(x),\\
 A_{n}^2= \begin{cases}
-\dfrac{1}{4}\dfrac{(n-j-a)(n-j-b)}{(n-j)}, & n \text { even}, \\[0.2cm]
-\dfrac{1}{4}(n+1), & n \text { odd},\ n\neq j, \\[0.2cm]
-\dfrac{1}{2}(1-\alpha)(j+1), & n \text { odd},\ n= j,
\end{cases}\\
C_{n}^2= \begin{cases}
\dfrac{1}{4}\dfrac{(n-j+a)(n-j+b)}{(n-j)}, & n \text { even}, \\[0.2cm]
\dfrac{1}{4}(n-2j-1), & n \text { odd},\ n\neq j, \\[0.2cm]
-\dfrac{1}{2}\alpha (j+1), & n \text { odd},\ n= j.
\end{cases}
\end{gather*}
\textit{Positivity condition}
\begin{gather*}
 a\leq -j, \qquad
 |b|\leq 1 , \qquad
 0 \leq \alpha \leq 1 ,
\qquad \text { or } \qquad
 |a|\leq 1 ,\qquad
b\leq -j, \qquad
0 \leq \alpha \leq 1.
\end{gather*}
\textit{Dunkl-difference equation}
\begin{gather*}
 \mathcal{D}_{\beta} \mathbf{P}_{n}^{(2)}(x)=\Lambda_{n}^{(\beta)} \mathbf{P}_{n}^{(2)}(x),\\
 \Lambda_{2n}^{(\beta)} = n(n-j),\qquad \Lambda_{2n+1}^{(\beta)} = n(n+1-j)+\beta,\\
 \mathcal{D}_{\beta} = \mathcal{D}_{0}+\beta \frac{(4x+b+a+j+1)}{8x}(I-R),\\
 \mathcal{D}_{0}=A(x) T^{+}+B(x) T^{-}+C(x) R+D(x) T^{+} R-(A(x)+B(x)+C(x)+D(x))I,
\end{gather*}
where $T^{\pm}f(x)=f(x\pm1)$ and $Rf(x)=f(-x)$ and where
\begin{gather*}
A(x)=\frac{ (x+\rho_{1}+1 )\bigl(x-\rho_{1}+\frac{1-j}{2}\bigl) (2 x-2 r_{1}+1 ) (2 x+2 r_{1}-j )}{8(x+1)(2 x+1)}, \\
B(x)=\frac{\bigl(x+\rho_{1}+\frac{j+1}{2}\bigr) (x-\rho_{1}-1 ) (2 x+2 r_{1}-1 ) (2 x-2 r_{1}+j )}{8 x(2 x-1)}, \\
C(x)=\frac{\bigl(x+\rho_{1}+\frac{j+1}{2}\bigr)\bigl(4 x^{2}+\omega\bigr)}{8 x}\\
\phantom{C(x)=}{}-\frac{\bigl(x+\rho_{1}+\frac{j+1}{2}\bigr) (x+\rho_{1}+1 ) (2 x-2 r_{1}+1 ) (2 x+2 r_{1}-j )}{8 x(2 x+1)}-B(x), \\
D(x)=-\frac{(2\rho_{1}+j+1) (x+\rho_{1}+1 ) (2 x-2 r_{1}+1 ) (2 x+2 r_{1}-j )}{16 x(x+1)(2 x+1)},\\
\omega = 4(1+\rho_1-r_1)\pp{\frac{1-j}{2}}+4r_1(1-r_1)-j,
\end{gather*}
where the change of parameter \eqref{CA_change_param} must be used.

\textit{Explicit expression}
\begin{gather*}
 \mathbf{P}_{2n}^{(2)}(x) =\kappa_n^{(1)}\sum_{k=0}^{n} A_{n,k},\qquad \mathbf{P}_{2n+1}^{(2)}(x) = \kappa_n^{(2)}\pp{x+\frac{b+j+1+a}{4}}\sum_{k=0}^{n} B_{n,k},
\end{gather*}
where
\begin{gather*}%\label{CA_Ank} \label{CA_kap1}
 A_{n,k}= \begin{cases}
\dfrac{(-n)_{k}(n-j)_{k}\bigl(-\frac{b+j+1+a}{4}+x\bigr)_{k}\bigl(-\frac{b+j+1+a}{4}-x\bigr)_{k}}{(1)_{k} \bigl(-\frac{b+j}{2}\bigr)_{k}\bigl(-\frac{j+a}{2}\bigr)_{k}\bigl(\frac{1-j}{2}\bigr)_{k}},\\
\hspace*{100mm} k\leq \frac{(j-1)}{2},\ k\leq n ,\\[0.3cm]
\pp{\dfrac{1}{\alpha}}\dfrac{(-n)_{k}\bigl(-\frac{b+j+1+a}{4}+x\bigr)_{k}\bigl(-\frac{b+j+1+a}{4}-x\bigr)_{k}}{(1)_{k} \bigl(-\frac{b+j}{2}\bigr)_{k}\bigl(-\frac{j+a}{2}\bigr)_{k}}\frac{(n-j)_{j-n}(1)_{n+k-j-1}}{\bigl(\frac{1-j}{2}\bigr)_{\frac{j-1}{2}}(1)_{k-(j+1)/2}}, \\
\hspace*{100mm} k> \frac{(j-1)}{2}, \ k\leq n,
\end{cases}\\
 \kappa_{n}^{(1)}= \begin{cases}
\dfrac{\bigl(-\frac{b+j}{2}\bigr)_{n}\bigl(-\frac{j+a}{2}\bigr)_{n}\bigl(\frac{1-j}{2}\bigr)_{n}}{(n-j)_{n}}, & n\leq \frac{(j-1)}{2}, \\[0.3cm]
(\alpha)\dfrac{\bigl(-\frac{b+j}{2}\bigr)_{n}\bigl(-\frac{j+a}{2}\bigr)_{n}\bigl(\frac{1-j}{2}\bigr)_{\frac{j-1}{2}}(1)_{n-(j+1)/2}}{(n-j)_{j-n}(1)_{2n-j-1}}, & n> \frac{(j-1)}{2},
\end{cases}
\end{gather*}
and
\begin{gather*}
 B_{n,k}= \begin{cases}
\dfrac{(-n)_{k}(n+1-j)_{k}\bigl(-\frac{b+j-3+a}{4}+x\bigr)_{k}\bigl(-\frac{b+j-3+a}{4}-x\bigr)_{k}}{(1)_{k}\bigl(\frac{3-j}{2}\bigr)_{k} \bigl(\frac{2-j-a}{2}\bigr)_{k}\bigl(\frac{2-j-b}{2}\bigr)_{k}}, & k< \frac{j-1}{2}, \\[0.3cm]
\pp{\dfrac{1}{\alpha}}\dfrac{(-n)_{k}\bigl(-\frac{b+j-3+a}{4}+x\bigr)_{k}\bigl(-\frac{b+j-3+a}{4}-x\bigr)_{k}}{(1)_{k}\bigl(\frac{2-j-a}{2}\bigr)_{k} \bigl(\frac{2-j-b}{2}\bigr)_{k}}\frac{(n+1-j)_{j-n-1}(1)_{n+k-j}}{\bigl(\frac{3-j}{2}\bigr)_{\frac{j-3}{2}}(1)_{k-(j-1)/2}}, & k\geq \frac{j-1}{2},
\end{cases} %\label{CA_Bnk}
\\
 \kappa_{n}^{(2)}= \begin{cases}
\dfrac{\bigl(\frac{2-j-a}{2}\bigr)_{n}\bigl(\frac{2-j-b}{2}\bigr)_{n}\bigl(\frac{3-j}{2}\bigr)_{n}} {(n+1-j)_{n}}, & n< \frac{j-1}{2}, \\[0.3cm]
(\alpha)\dfrac{\bigl(\frac{2-j-a}{2}\bigr)_{n}\bigl(\frac{2-j-b}{2}\bigr)_{n}\bigl(\frac{3-j}{2}\bigr)_{\frac{j-3}{2}}(1)_{n-(j-1)/2}}{(n+1-j)_{j-n-1}(1)_{2n-j}}, & n\geq \frac{j-1}{2}.
\end{cases}%\label{CA_kap2}
\end{gather*}
\textit{Orthogonality relation}
\begin{gather*}
 x_{4s} = -\pp{\frac{4s-j+a-b}{4}}-\frac{1}{4}, \qquad s\in \{0,1,\dots,(j-1)/2\},\\
 x_{4s+1} = \pp{\frac{4s-j-a-b}{4}}-\frac{1}{4},\qquad s\in \{0,1,\dots,(j-1)/2\},\\
 x_{4s+2} = \pp{\frac{4s+2-j+a-b}{4}}-\frac{1}{4},\qquad s\in \{0,1,\dots,(j-1)/2\},\\
 x_{4s+3} = -\pp{\frac{4s+2-j-a-b}{4}}-\frac{1}{4},\qquad s\in \{0,1,\dots,(j-3)/2\},\\
 w_{4s} = \frac{2\alpha h_{2j}}{(-s)_{s}(1)_{(j-1)/2-s}\bigl(-s+\frac{1-a}{2}\bigr)_{(j-1)/2}\bigl(s+\frac{1+a-b-j}{2}\bigr)_{(j+1)/2}\bigl(s-\frac{b+j}{2}\bigr)_{(j+1)/2}},\\
 w_{4s+1} = \frac{2(1-\alpha)h_{2j}}{(-s)_{s}(1)_{(j-1)/2-s}\bigl(-s+\frac{1+a}{2}\bigr)_{(j+1)/2}\bigl(s+\frac{1-a-b-j}{2}\bigr)_{(j-1)/2}\bigl(s-\frac{b+j}{2}\bigr)_{(j+1)/2}},\\
 w_{4s+2} = -\frac{2\alpha h_{2j}}{(-s)_{s}(1)_{(j-1)/2-s}\bigl(-s-\frac{1+a}{2}\bigr)_{(j+1)/2}\bigl(s+\frac{1+a-b-j}{2}\bigr)_{(j+1)/2}\bigl(s+\frac{2-b-j}{2}\bigr)_{(j-1)/2}},\\
 w_{4s+3} = -\frac{2(1-\alpha)h_{2j}}{(-s)_{s}(1)_{(j-3)/2-s}\bigl(-s-\frac{1-a}{2}\bigr)_{(j+1)/2}\bigl(s+\frac{1-a-b-j}{2}\bigr)_{(j+1)/2}\bigl(s+\frac{2-b-j}{2}\bigr)_{(j+1)/2}},\\
 h_{2n} = \frac{\sqrt{(-j)_{n}(1)_{n}\bigl(\frac{2+a-j}{2}\bigr)_{n}\bigl(-\frac{j+a}{2}\bigr)_{n} \bigl(\frac{2+b-j}{2}\bigr)_{n}\bigl(-\frac{j+b}{2}\bigr)_{n}}}{2^{2n}\sqrt{\bigl(\frac{2-j}{2}\bigr)_{n}\bigl(-\frac{j}{2}\bigr)_{n}}}
 \times \begin{cases}
1, & n\leq \frac{(j-1)}{2}, \\[0.1cm]
2\sqrt{\alpha(1-\alpha)}, & n> \frac{(j-1)}{2},
\end{cases}\\
h_{2n+1} = \frac{\sqrt{(-1)(-j)_{n+1}(1)_{n}\bigl(\frac{2+a-j}{2}\bigr)_{n}\bigl(-\frac{j+a}{2}\bigr)_{n+1}\bigl(\frac{2+b-j}{2}\bigr)_{n} \bigl(-\frac{j+b}{2}\bigr)_{n+1}}}{2^{2n+1}\sqrt{\bigl(\frac{2-j}{2}\bigr)_{n}\bigl(-\frac{j}{2}\bigr)_{n+1}}}\\ \phantom{h_{2n+1} =}{}\times \begin{cases}
1, & n<\frac{(j-1)}{2},\\[0.1cm]
\sqrt{2\alpha}, & n=\frac{(j-1)}{2},\\[0.1cm]
2\sqrt{\alpha(1-\alpha)}, & n>\frac{(j-1)}{2}.
\end{cases}
\end{gather*}
The orthogonality relation reads
\begin{gather*}
 \sum_{s=0}^{2j}w_{s} \mathbf{P}_{n}^{(2)}(x_s)\mathbf{P}_{m}^{(2)}(x_s) = h_n^2 \delta_{nm}\label{CA_ortho}.
\end{gather*}

\begin{figure}[t]\centering
\begin{tikzpicture}[/pgf/declare function={p = 1.5;},/pgf/declare function={q = 2;}]
\draw[-,ultra thick] (-7.5,0)--(-5.5,0);
\draw[dashed,ultra thick] (-5.75,0)--(-4.25,0);
\draw[-,ultra thick] (-4.25,0)--(4.25,0);
\draw[dashed,ultra thick] (4.25,0)--(5.75,0);
\draw[->,ultra thick] (5.5,0)--(7.5,0)node[right]{$x$};
%\draw[->,ultra thick] (-8,0)--(8,0) node[right]{$x$};
\draw[-,ultra thick] (0,-0.1) node[below]{$0$} --(0,0.1);
\foreach \Point in {(-p*9/8,0),(-p*11/8,0),(-p*17/8,0),(-p*19/8,0),(-q*25/8,0),(-q*27/8,0),(p*3/8,0),(p*9/8,0),(p*11/8,0),(p*17/8,0),(p*19/8,0),(q*25/8,0),(q*27/8,0)}{
 \node[red] at \Point {\textbullet};}
\draw[<->] (-p*9/8,0.5) -- (-0.05,0.5) node [midway, above] (1) {$d_1$};
\draw[<->] (-p*17/8,0.5) -- (-p*11/8,0.5) node [midway, above] (11) {$d_4$};
\draw[<->] (-p*9/8,-0.5) -- (-p*11/8,-0.5) node [midway, below] (7) {$d_3$};
\draw[<->] (-p*19/8,-0.5) -- (-p*17/8,-0.5) node [midway, below] (8) {$d_3$};
\draw[<->] (-q*25/8,-0.5) -- (-q*27/8,-0.5) node [midway, below] (6) {$d_3$};
\draw[<->] (0.05,0.5) -- (p*3/8,0.5) node [midway, above] (2) {$d_2$};
\draw[<->] (p*3/8,-0.5) -- (p*9/8,-0.5) node [midway, below] (9) {$d_4$};
\draw[<->] (p*9/8,0.5) -- (p*11/8,0.5) node [midway, above] (3) {$d_3$};
\draw[<->] (p*17/8,0.5) -- (p*19/8,0.5) node [midway, above] (4) {$d_3$};
\draw[<->] (p*17/8,-0.5) -- (p*11/8,-0.5) node [midway, below] (10) {$d_4$};
\draw[<->] (q*25/8,0.5) -- (q*27/8,0.5) node [midway, above] (5) {$d_3$};
\end{tikzpicture}
\caption{Depiction of a Bannai--Ito bi-lattice.}
\end{figure}
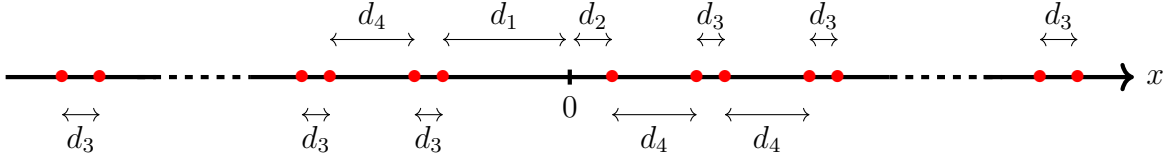

\begin{figure}[t]\centering
\begin{tikzpicture}[/pgf/declare function={p = 1.5;},/pgf/declare function={q = 2;}]
\draw[-,ultra thick] (-7.5,0)--(-5.5,0);
\draw[dashed,ultra thick] (-5.75,0)--(-4.25,0);
\draw[-,ultra thick] (-4.25,0)--(4.25,0);
\draw[dashed,ultra thick] (4.25,0)--(5.75,0);
\draw[->,ultra thick] (5.5,0)--(7.5,0)node[right]{$x$};
%\draw[->,ultra thick] (-8,0)--(8,0) node[right]{$x$};
\draw[-,ultra thick] (0,-0.1) node[below]{$0$} --(0,0.1);
\foreach \Point in {(-p*8/8,0),(-p*12/8,0),(-p*16/8,0),(-p*20/8,0),(-p*24/8-p*8/8,0),(-p*24/8-p*12/8,0),(p*4/8,0),(p*8/8,0),(p*12/8,0),(p*16/8,0),(p*20/8,0),(p*24/8+p*8/8,0),(p*24/8+p*12/8,0)}{
 \node[red] at \Point {\textbullet};}
\draw[<->] (-p*8/8,0.5) -- (-0.05,0.5) node [midway, above] (1) {$d_1$};
\draw[<->] (-p*16/8,0.5) -- (-p*12/8,0.5) node [midway, above] (11) {$d_3$};
\draw[<->] (-p*8/8,-0.5) -- (-p*12/8,-0.5) node [midway, below] (7) {$d_3$};
\draw[<->] (-p*20/8,-0.5) -- (-p*16/8,-0.5) node [midway, below] (8) {$d_3$};
\draw[<->] (-p*24/8-p*8/8,-0.5) -- (-p*24/8-p*12/8,-0.5) node [midway, below] (6) {$d_3$};
\draw[<->] (0.05,0.5) -- (p*4/8,0.5) node [midway, above] (2) {$d_2$};
\draw[<->] (p*4/8,-0.5) -- (p*8/8,-0.5) node [midway, below] (9) {$d_3$};
\draw[<->] (p*8/8,0.5) -- (p*12/8,0.5) node [midway, above] (3) {$d_3$};
\draw[<->] (p*16/8,0.5) -- (p*20/8,0.5) node [midway, above] (4) {$d_3$};
\draw[<->] (p*16/8,-0.5) -- (p*12/8,-0.5) node [midway, below] (10) {$d_3$};
\draw[<->] (p*24/8+p*8/8,0.5) -- (p*24/8+p*12/8,0.5) node [midway, above] (5) {$d_3$};
\end{tikzpicture}
\caption{Reduction of the Bannai--Ito bi-lattice to a single Bannai--Ito lattice.}
\end{figure}
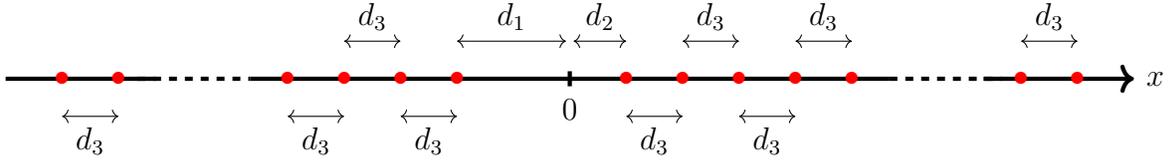

\begin{figure}[t!]\centering
\begin{tikzpicture}[/pgf/declare function={p = 1.5;},/pgf/declare function={q = 2;}]
\draw[-,ultra thick] (-7.5,0)--(-5.5,0);
\draw[dashed,ultra thick] (-5.75,0)--(-4.25,0);
\draw[-,ultra thick] (-4.25,0)--(4.25,0);
\draw[dashed,ultra thick] (4.25,0)--(5.75,0);
\draw[->,ultra thick] (5.5,0)--(7.5,0)node[right]{$x$};
%\draw[->,ultra thick] (-8,0)--(8,0) node[right]{$x$};
\draw[-,ultra thick] (0,-0.1) node[below]{$0$} --(0,0.1);
\foreach \Point in {(-p*4/8,0),(-p*8/8,0),(-p*12/8,0),(-p*16/8,0),(-p*20/8,0),(-p*24/8-p*8/8,0),(-p*24/8-p*12/8,0),(0,0),(p*4/8,0),(p*8/8,0),(p*12/8,0),(p*16/8,0),(p*20/8,0),(p*24/8+p*8/8,0),(p*24/8+p*12/8,0)}{
 \node[red] at \Point {\textbullet};}
\draw[<->] (-p*4/8,0.5) -- (-0.05,0.5) node [midway, above] (1) {$d$};
\draw[<->] (-p*8/8,-0.5) -- (-p*4/8,-0.5) node [midway, below] (11) {$d$};
\draw[<->] (-p*8/8,+0.5) -- (-p*12/8,+0.5) node [midway, above] (7) {$d$};
\draw[<->] (-p*24/8-p*8/8,-0.5) -- (-p*24/8-p*12/8,-0.5) node [midway, below] (6) {$d$};
\draw[<->] (0.05,0.5) -- (p*4/8,0.5) node [midway, above] (2) {$d$};
\draw[<->] (p*4/8,-0.5) -- (p*8/8,-0.5) node [midway, below] (9) {$d$};
\draw[<->] (p*8/8,0.5) -- (p*12/8,0.5) node [midway, above] (3) {$d$};
\draw[<->] (p*24/8+p*8/8,0.5) -- (p*24/8+p*12/8,0.5) node [midway, above] (5) {$d$};
\end{tikzpicture}
\caption{Reduction of the single Bannai--Ito lattice into a single linear lattice.}
\end{figure}
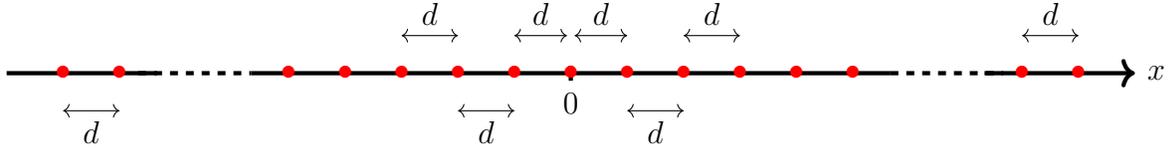

\section[Para-Bannai--Ito for N=2j+1, j odd]{Para-Bannai--Ito for $\boldsymbol{N=2j+1}$, $\boldsymbol{j}$ odd}\label{sec:appB}
\textit{Truncation from the general Bannai--Ito polynomials}
\begin{gather*}
 \rho_1+\rho_2 = -\frac{j+1}{2}+e_1t,\qquad
 -r_1-r_2 = -\frac{j+1}{2}+e_2t\label{CB_param_jodd}.
\end{gather*}
The para-Bannai--Ito polynomials are obtained through a $t\to0$ limit, with
\begin{gather*}%\label{CB_change_param}
 \rho_1 = \frac{b-j-1+a}{4},\qquad -r_1 = \frac{b-j-1-a}{4},\qquad \frac{e_1}{e_1+e_2} = \alpha ,\qquad \frac{e_2}{e_1+e_2} = (1-\alpha).
\end{gather*}
\textit{Recurrence relation}
\begin{gather*}%\label{CB_RRPBI_2j}\label{CB_RRPBICo_2j}
 x\mathbf{P}_{n}^{(3)}(x)=\mathbf{P}_{n+1}^{(3)}(x)+\pp{\frac{b-j-1+a}{4}-A_n^3-C_n^3}\mathbf{P}_{n}^{(3)}(x)+A_{n-1}^3C_n^3\mathbf{P}_{n-1}^{(3)}(x),\\
 A_{n}^3=\begin{cases}
\dfrac{1}{4}\dfrac{(n-j+a)(n-j+b)}{(n-j)}, & n \text { even}, \\[0.2cm]
\dfrac{1}{4}(n-2j-1), & n \text { odd},\ n\neq j, \\[0.2cm]
-\dfrac{1}{2}\alpha (j+1), & n \text { odd},\ n= j,
\end{cases}\\
C_{n}^3=\begin{cases}
-\dfrac{1}{4}n, & n \text { even},\ n\neq j+1, \\[0.2cm]
-\dfrac{1}{4}\dfrac{(n-j-1-a)(n-j-1-b)}{(n-j-1)}, & n \text { odd}, \\[0.2cm]
-\dfrac{1}{2}(1-\alpha)(j+1), & n \text { even},\ n= j+1.
\end{cases}
\end{gather*}
\textit{Positivity condition}
\begin{gather*}
 |a|\geq j, \qquad
 |b|\leq 1 ,\qquad
0 \leq \alpha \leq 1 , \qquad \text { or } \qquad
 |a|\leq 1 ,\qquad
|b|\geq j, \qquad
0 \leq \alpha \leq 1.
\end{gather*}
\textit{Dunkl-difference equation}
\begin{gather*}
L \mathbf{P}_{n}^{(1)}(x)=\lambda_{n} \mathbf{P}_{n}^{(1)}(x),\qquad
 \lambda_{2n} = n ,\qquad \lambda_{2n+1} = (j-n),\\
 L=F(x)(I-R)+G(x)\left(T^{+} R-I\right),
\end{gather*}
where
\begin{gather*}
 G(x) = \frac{(4x+1+b-a-j)(4x+1+a-b-j)}{16(2x+1)},\\ F(x) = \frac{(4x+j+1-a-b)(4x+j+1+a+b)}{32x}.
\end{gather*}
\textit{Explicit expression}
\begin{gather*}
 \mathbf{P}_{n}^{(3)}(x) = \mathbf{P}_{n}^{(2)}(x)-C_n^{3}\mathbf{P}_{n-1}^{(2)}(x).
\end{gather*}
\textit{Orthogonality relation}
\begin{gather*}
 x_{4s} = -\pp{\frac{4s-j+a-b}{4}}-\frac{1}{4}, \qquad s\in \{0,1,\dots,(j-1)/2\},\\
 x_{4s+1} = \pp{\frac{4s-j-a-b}{4}}-\frac{1}{4},\qquad s\in \{0,1,\dots,(j-1)/2\},\\
 x_{4s+2} = \pp{\frac{4s+2-j+a-b}{4}}-\frac{1}{4},\qquad s\in \{0,1,\dots,(j-1)/2\},\\
 x_{4s+3} = -\pp{\frac{4s+2-j-a-b}{4}}-\frac{1}{4},\qquad s\in \{0,1,\dots,(j-1)/2\},\\
 w_{4s} = \frac{2\alpha h_{2j}}{(-s)_{s}(1)_{(j-1)/2-s}\bigl(-s+\frac{1-a}{2}\bigr)_{(j+1)/2}\bigl(s+\frac{1+a-b-j}{2}\bigr)_{(j+1)/2}\bigl(s-\frac{b+j}{2}\bigr)_{(j+1)/2}},\\
 w_{4s+1} = \frac{2(1-\alpha)h_{2j}}{(-s)_{s}(1)_{(j-1)/2-s}\bigl(-s+\frac{1+a}{2}\bigr)_{(j+1)/2}\bigl(s+\frac{1-a-b-j}{2}\bigr)_{(j+1)/2}\bigl(s-\frac{b+j}{2}\bigr)_{(j+1)/2}},\\
 w_{4s+2} = -\frac{2\alpha h_{2j}}{(-s)_{s}(1)_{(j-1)/2-s}\bigl(-s-\frac{1+a}{2}\bigr)_{(j+1)/2}\bigl(s+\frac{1+a-b-j}{2}\bigr)_{(j+1)/2}\bigl(s+\frac{2-b-j}{2}\bigr)_{(j+1)/2}},\\
 w_{4s+3} = -\frac{2(1-\alpha)h_{2j}}{(-s)_{s}(1)_{(j-1)/2-s}\bigl(-s-\frac{1-a}{2}\bigr)_{(j+1)/2}\bigl(s+\frac{1-a-b-j}{2}\bigr)_{(j+1)/2}\bigl(s+\frac{2-b-j}{2}\bigr)_{(j+1)/2}},\\
 h_{2n} = \frac{\sqrt{(-j)_{n}(1)_{n}\bigl(-\frac{j-a}{2}\bigr)_{n}\bigl(-\frac{j+a}{2}\bigr)_{n}\bigl(-\frac{j-b}{2}\bigr)_{n} \bigl(-\frac{j+b}{2}\bigr)_{n}}}{2^{2n}\bigl(-\frac{j}{2}\bigr)_{n}} \times \begin{cases}
1, & n< \frac{(j+1)}{2}, \\[0.1cm]
2\sqrt{\alpha(1-\alpha)}, & n\geq \frac{(j+1)}{2},
\end{cases}\\
 h_{2n+1} = \frac{\sqrt{(-1)(-j)_{n}(1)_{n}\bigl(-\frac{j-a}{2}\bigr)_{n+1}\bigl(-\frac{j+a}{2}\bigr)_{n+1}\bigl(-\frac{j-b}{2}\bigr)_{n+1} \bigl(-\frac{j+b}{2}\bigr)_{n+1}}}{2^{2n+1}\bigl(-\frac{j}{2}\bigr)_{n+1}}\\ \qquad \qquad\times \begin{cases}
1, & n< \frac{(j+1)}{2}, \\[0.1cm]
2\sqrt{\alpha(1-\alpha)}, & n\geq \frac{(j+1)}{2}.
\end{cases}
\end{gather*}
The orthogonality relation reads
\begin{gather*}
 \sum_{s=0}^{2j+1}w_{s} \mathbf{P}_{n}^{(3)}(x_s)\mathbf{P}_{m}^{(3)}(x_s) = h_n^2 \delta_{nm}\label{CB_ortho}.
\end{gather*}

\section{Depiction of the Bannai--Ito bi-lattice}\label{secappC}
We show here an example of the Bannai--Ito bi-lattice for $a=-(c+j+1)$ and $c>0$, which correspond to the first set of restrictions of \eqref{eq:cond_posi} with $j$ even. We see that the parameter $c$, and consequently $a$, controls the gap between the positive and negative part of the grid. As for the parameter $b$, it controls the gap between consecutive points on each separate side. The first and last points are at the coordinate $x=\pm\frac{2j+b+c}{4}$ and the distances $d_3$ and $d_4$ alternate until they reach the end on both side.
We have precisely,
\begin{gather*}
 d_1 = \frac{c+2-b}{4},\qquad d_2 = \frac{b+c}{4},\qquad d_3 = \frac{1+b}{2},\qquad d_4=\frac{1-b}{2}.
\end{gather*}
When we use $b=0$, we obtain equally distanced points on each side. This is equivalent to a~single Bannai--Ito lattice (or two separate linear grid).

Now, we have
\begin{gather*}
 d_1 = \frac{c}{4}+\frac{1}{2},\qquad d_2 = \frac{c}{4},\qquad d_3 = \frac{1}{2}.
\end{gather*}
Finally, if we use $c=0$ (or $a=-j-1$), we obtain one simple linear grid where the points are all separated by $d=1/2$.

\subsubsection*{Acknowledgements}
JP holds a scholarship from Fonds de recherche qu\'eb\'ecois -- nature et technologies (FRQNT) and an academic excellence scholarship from Hydro-Qu\'ebec. The research of LV is supported in part by a Discovery grant from the Natural Science and Engineering Research Council (NSERC) of Canada.

\pdfbookmark[1]{References}{ref}
\LastPageEnding

\end{document}